\documentclass[11pt]{article}
\usepackage[english]{style}
\usepackage{latexsym,amsmath,amssymb,amsfonts,diagrams,epsf,fancyvrb}

\DeclareFontFamily{U}{rsf}{}
\DeclareFontShape{U}{rsf}{m}{n}{
  <5> <6> rsfs5 <7> <8> <9> rsfs7 <10-> rsfs10}{}
\DeclareMathAlphabet{\mathscr}{U}{rsf}{m}{n}

\DeclareMathAlphabet{\mathgth}{U}{euf}{m}{n}

\DeclareFontFamily{U}{cyr}{}
\DeclareFontShape{U}{cyr}{m}{n}{
  <5> wncyr5 <6> wncyr6 <7> wncyr7 <8> wncyr8 <9> wncyr9 <10-> wncyr10}{}
\DeclareMathAlphabet{\mathcyr}{U}{cyr}{m}{n}

\input cyracc.def

\makeatletter
\def\operator@font{\sf}
\makeatother

\newcommand{\cA}{{\mathscr A}}

\newcommand{\cC}{{\mathscr C}}
\newcommand{\cE}{{\mathscr E}}

\newcommand{\cI}{{\mathscr I}}

\newcommand{\cN}{{\mathscr N}}
\newcommand{\cO}{{\mathscr O}}

\newcommand{\tS}{\tilde{S}}

\newcommand{\D}{{\mathbf D}_{\mathsf{coh}}^b}
\newcommand{\Pf}{\mathbf{Pf}}

\newcommand{\Fuk}{{\mathbf{Fuk}}}

\newcommand{\R}{\mathbf{R}}
\newcommand{\Ld}{\mathbf{L}}
\newcommand{\lotimes}{\stackrel{\Ld}{\otimes}}

\newcommand{\smooth}{{\mathsf{sm}}}

\DeclareMathOperator{\Cone}{Cone}
\DeclareMathOperator{\Lie}{Lie}

\DeclareMathOperator{\Ann}{Ann}
\DeclareMathOperator{\Sym}{Sym}

\newcommand{\Hom}{{\mathsf{Hom}}}
\DeclareMathOperator{\End}{End}
\DeclareMathOperator{\Koszul}{Koszul}
\DeclareMathOperator{\Tor}{Tor}
\DeclareMathOperator{\Ker}{Ker}

\DeclareMathOperator{\Pic}{Pic}

\DeclareMathOperator{\GL}{GL}

\DeclareMathOperator{\rk}{rk}
\DeclareMathOperator{\Ext}{Ext}

\DeclareMathOperator{\Span}{Span}

\newcommand{\T}{\mathbf{T}}

\newcommand{\ra}{\rightarrow}

\newcommand{\lra}{\longrightarrow}

\newcommand{\C}{\mathbf{C}}
\newcommand{\G}{\mathbf{G}}

\newcommand{\Z}{\mathbf{Z}}

\newcommand{\iso}{\cong}

\newcommand{\pj}{\mathbf{P}}
\newcommand{\PW}{\mathbf{W}}
\newcommand{\PM}{\mathbf{M}}

\newarrow{Equal} =====

\renewcommand{\phi}{\varphi}

\author{%
Lev Borisov, Andrei C\u ald\u araru\thanks{Mathematics Department,
University of Wisconsin--Madison, 480 Lincoln Drive, Madison, WI
53706--1388, USA, {\em e-mail: }{\tt andreic@math.wisc.edu,
borisov@math.wisc.edu}}}

\title{The Pfaffian-Grassmannian derived equivalence}

\date{}

\begin{document}

\maketitle

\begin{abstract}
  We argue that there exists a derived equivalence between Calabi-Yau
  threefolds obtained by taking dual linear sections (of the
  appropriate codimension) of the Grassmannian $\G(2,7)$ and the
  Pfaffian $\Pf(7)$.  The existence of such an equivalence has been
  conjectured by physicists for almost ten years, as the two families
  of Calabi-Yau threefolds are believed to have the same mirror.  It
  is the first example of a derived equivalence between Calabi-Yau
  threefolds which are provably non-birational.
\end{abstract}

\section*{Introduction}

\paragraph
Let $V$ be a vector space of dimension seven over $\C$ (or any
algebraically closed field of characteristic zero), and let
\[ \G = \G(2, V) \] 
be the Grassmannian of planes in $V$.  The Pl\"ucker map embeds $\G$
as a smooth subvariety of dimension $10$ of 
\[ \pj = \pj^{20} = \pj(\wedge^2 V). \]
Regard the dual projective space 
\[ \pj^* = \pj(\wedge^2 V^*) \] 
as the projectivization of the space of two-forms on $V$.  The
Pfaffian locus 
\[ \Pf \subset \pj^*\] 
is defined to be the projectivization of the locus of degenerate
two-forms on $V$ (forms of rank $\leq 4$).  Equations for $\Pf$ can be
obtained by taking the Pfaffians of the diagonal minors of a
skew-symmetric $7 \times 7$ matrix of linear forms on $V$.

While the Grassmannian $\G$ is smooth, the Pfaffian $\Pf$ is a
singular subvariety of $\pj^*$ of dimension $17$.  Indeed, a point
$\omega\in \Pf$ will be singular precisely when the rank of $\omega$
is two.  (Recall that the rank of a two-form $\omega$ -- or,
equivalently, of a skew symmetric matrix -- is always even.  In our
case, since we eliminate $\omega=0$ by projectivizing, this rank could
be two, four, or six.  The general two-form on $V$ has rank six; it has
rank four at the smooth points of $\Pf$, and it has rank two at its
singular points, which coincide with $\G(2, V^*)$ in its Pl\"ucker
embedding.)

The Pfaffian is the classical projective dual of the Grassmannian: 
\[ \Pf = \{ y \in \pj^*~:~ \G \cap H_y\mbox{ is singular} \}, \]
where $H_y$ is the linear space in $\pj$ corresponding to $y$.

\paragraph
Consider a seven-dimensional linear subspace
\[ W\subset \wedge^2 V^*, \] 
and denote by $\PW$ its image in $\pj^*$.  Let $Y$ be the
intersection of $\PW$ with $\Pf$.

On the dual side, let 
\[ M = \Ann(W) \subset \wedge^2 V \] 
be the $14$-dimensional annihilator of $W$; again, we will use bold-face
$\PM$ to denote its projectivization in
$\pj$, which has codimension seven.  Let $X$ be the intersection of
$\PM$ and $\G$.  

\smallskip\noindent
The main result of this paper is the following theorem.

\begin{Theorem}
\label{thm:mainthm}
  For a given choice of $W$, if either $X$ or $Y$ has dimension three,
  then $X$ is smooth if and only if $Y$ is.  When this happens, $X$
  and $Y$ are Calabi-Yau threefolds with
  \[ h^{1,1} = 1,\quad h^{1,2} = 50, \] 
  and there exists an equivalence of derived categories
  \[ \Phi : \D(Y) \stackrel{\sim}{\lra} \D(X). \]
\end{Theorem}

\paragraph
Such a result has been conjectured for a while.  Indeed, R\o
dland~\cite{Rod} argued, by comparing solutions to the Picard-Fuchs
equation, that the families of $X$'s and of $Y$'s appear to have the
same mirror family.  Recently, Hori and Tong~\cite{HorTon} gave a more
detailed string theory argument supporting the same conclusion.  If we
denote by $Z$ their common mirror, Kontsevich's Homological Mirror
Symmetry conjecture predicts
\[ \D(X) \iso \Fuk(Z) \iso \D(Y). \]

\paragraph
This appears to be the first example of a derived equivalence between
Calabi-Yau threefolds which can be proved to be non-birational.
Indeed, if $X$ and $Y$ were birational, they would have to differ by a
sequence of flops because they are minimal in the sense of Mori
theory.  On the other hand, no flops are possible on either $X$ or
$Y$, because they have Picard number $h^{1,1} = 1$.  Therefore they
would have to be isomorphic.  However, this can not be true, since if
we denote by $H_X$ and $H_Y$, respectively, the ample generators of
the Picard groups of $X$ and of $Y$, we have~\cite{Rod}
\[ H_X^3 = 42, \qquad H_Y^3 = 14. \]

\paragraph
Our results appear to fit very well with the theory of Homological
Projective Duality developed by Kuznetsov~\cite{Kuz}.  Indeed, we have
a pair of varieties $\G$ and $\Pf$, embedded in dual projective
spaces, whose dual linear sections are derived equivalent.  It
would be interesting to understand this relationship further.  To this
end, we make some comments in Section~\ref{sec:comments} about how this
example seems to fit in the general theory.

\paragraph
Let us now briefly describe the construction of the derived
equivalence $\Phi$.  Recall~\cite{BriFM} that, in order to give $\Phi$,
it is essentially enough to describe the image of $\Phi$ on structure
sheaves of points $\cO_y$ for $y \in Y$, and to check that the family
$\{\Phi\cO_y\}_{y\in Y}$ is an orthonormal basis for $\D(X)$.  

A point $y$ in $Y$ can be regarded as a two-form on $V$ of rank four.  As
such it has a kernel $K$ which is a three-dimensional linear subspace of
$V$.  On the other hand, a point $x$ in $X$ corresponds to a
two-dimensional linear subspace $T$ of $V$.  For a general choice of $x$
and $y$ we will have
\[ T \cap K = 0. \]
However, for a fixed $y\in Y$, the set of points $x\in X$ for
which $T$ intersects $K$ non-trivially is a curve $C_y$ in
$X$.  Our choice for the functor $\Phi$ is to set 
\[ \Phi\cO_y = \cI_{C_y}, \]
where $\cI_{C_y}$ denotes the ideal sheaf of the curve $C_y$.  The
technical core of the paper consists then of showing that this choice
gives rise to an orthonormal family.

It would be interesting to find out the relationship of this approach
to the construction of Donaldson-Thomas moduli spaces.  In the course
of our proof we see that $Y$ is such a moduli space of ideal
sheaves of curves on $X$.  Perhaps other examples of derived
equivalences can be obtained as Donaldson-Thomas spaces.

\paragraph
The paper is organized as follows.  In Section~\ref{sec:intro} we set
up the notation, and we recall the appropriate definitions from linear
algebra.  We also sketch the standard constructions of the tangent
spaces to the Pfaffian and Grassmannian varieties.  In
Section~\ref{sec:smoothness} we give a purely linear-algebra argument
for the fact that $X$ is smooth if and only if $Y$ is.  We study the
geometry of certain Schubert cycles in Section~\ref{sec:gvan}, and we
prove a vanishing result for $\Ext$ groups between them on the
Grassmannian.  In Section~\ref{sec:allcurves} we define the curves
$C_y$ which are parametrized by the points of $Y$, and we argue that
they all have dimension one.
Section~\ref{sec:redGrass} reduces the computation of the
orthogonality of the family $\{\cI_{C_y}\}$ to the vanishing result on
$\G$ proved earlier.  We define the functor $\Phi$ in
Section~\ref{sec:deq} and we argue that it is an equivalence.
Finally, in Section~\ref{sec:comments} we discuss connections with
results of Kuznetsov on Homological Projective Duality.

\paragraph {\bf Acknowledgments.} The second author would like to
thank Mark Gross for first introducing him to the physicists'
conjecture, to Kentaro Hori for pointing out the reference~\cite{Rod},
and to Ron Donagi for comments on the geometry of the Grassmannian.
Michael Stillman provided comments on the initial Macaulay
computations which confirmed our results.  Dan Abramovich had comments
on the relationship to Donaldson-Thomas invariants.  Alexander
Kuznetsov explained connections with Homological Projective Duality in
higher dimensions, and pointed out an alternative argument for
$\G$-vanishing.  This material is based upon work supported by the
National Science Foundation under Grants No.\ DMS-0456801 and
DMS-0556042.

\section{Basic facts from linear algebra}
\label{sec:intro}

In this section we set up the basic notation and review some standard
results from linear algebra and differential geometry.

\paragraph
Let $V$ be a vector space.  Denote by $\G(k, V)$ the Grassmannian of
$k$-dimensional planes in $V$.  The Pl\"ucker embedding is the map
\[ \G(k, V) \hookrightarrow \pj(\wedge^k V), \] 
which maps a point $x$ of $\G(k, V)$ (i.e., a $k$-plane $x\subset V$)
to the point
\[ [e_1 \wedge e_2 \wedge \cdots \wedge e_k] \in \pj(\wedge^k V), \]
where $e_1,\ldots,e_k$ is any basis of $x$. 

If $x$ is a point in the Grassmannian, the tangent space $T_{\G(k, V),
  x}$ is canonically identified with 
\[ \Hom(x, V/x). \]
Thus, if $V$ has dimension $n$, the Grassmannian $\G(k, V)$ has
dimension $k(n-k)$. 

In particular, if $k=1$, the tangent space $T_{\pj V, x}$ to the
projective space at a point $x$ is naturally identified with $\Hom(x,
V/x)$ which is (up to scalars) the same as $V/x$.

\paragraph
The tangent map to the Pl\"ucker embedding is the map
\[ \Hom(x, V/x) \ra \Hom(\wedge^k x, \wedge^k V/\wedge^k x), \]
given by
\[ \phi \mapsto \left (v_1\wedge\cdots \wedge v_k \mapsto \sum_{i=1}^k v_1
\wedge v_2 \wedge \cdots \wedge v_{i-1} \wedge \phi(v_i) \wedge
v_{i+1} \wedge \cdots \wedge v_k \right). \]

\paragraph
An element $y\in \wedge^k V^*$ will be called a $k$-form on $V$.
Since we have 
\[ \wedge^k V^* \iso \left ( \wedge^k V \right )^*, \] 
the annihilator of $y$, $\Ann(y)$, is a hyperplane $H\subset \wedge^k
V$.  It consists of all $x\in \wedge^k V$ such that $y(x) = 0$.

\paragraph
From now on we will concentrate on the case $k=2$, and define
\[ \pj = \pj(\wedge^2 V),\quad \pj^* = \pj(\wedge^2 V^*), \quad \G =
\G(2, V). \]
Frequently, when regarding a point $x\in \G$ as a two-plane in $V$ we
will call this plane $T$. 

If $y$ is a two-form on $V$, its kernel $K$ is defined to be the set of
all $v\in V$ such that $y(v\wedge w) = 0$ for all $w\in V$.  We
define the rank of $y$ by
\[ \rk y = n-\dim K. \]
The rank of a two-form is always even.

The kernel of a form $y$ does not change if we multiply $y$ by a
non-zero scalar.  Thus we can speak of the kernel of a point $y\in
\pj^*$, and we will frequently refer to such a point as a two-form.

\begin{Proposition}
\label{prop:hyptgrass}
  Let $y$ be a point in $\pj^*$, and let $H$ be the corresponding
  hyperplane in $\pj$.  Let $x$ be a point of intersection of $H$ and
  $\G$.  Then $H$ is tangent to $\G$ at $x$ (under the Pl\"ucker
  embedding) if and only if $T\subset K$, where $T$ is the two-plane in
  $V$ corresponding to $x$, and $K$ is the kernel of $y$, regarded as
  a two-form on $V$.
\end{Proposition}

\begin{Proof}
First assume $T\subset K$.  Let $\phi$ be a tangent vector to $\G$
at $x$, i.e., a morphism $\phi: T \ra V/T$.  Then the image of
$\phi$ under the differential of the Pl\"ucker map is the map
$\wedge^2 T \ra \wedge^2 V / \wedge^2 T$ given by
\[ t_1 \wedge t_2 \mapsto \phi(t_1) \wedge t_2 + t_1 \wedge
\phi(t_2). \]
Choosing $t_1$ and $t_2$ linearly independent vectors in $T$, we
identify this (up to scalars) with the element
\[  t = \phi(t_1) \wedge t_2 + t_1 \wedge \phi(t_2) \]
of $\wedge^2 V / \wedge^2 T$.  Since $T\subset K$, $t_1, t_2$ are in
the kernel of $y$, and therefore $y$ vanishes on $t$.  In other words,
we have proven that 
\[ \phi \in T_{\G, x} \implies \phi \in T_{H, x}, \]
or that $H$ is tangent to $\G$ at $x$. 

Conversely, assume $T\not\subset K$.  At least one of $t_1$, $t_2$ is
not in $K$, assume it is $t_1$.  There exists then $v\in V$ such that
$y(t_1 \wedge v) \neq 0$.  Define $\phi$ on the basis $t_1, t_2$ of
$T$ by setting
\begin{align*} 
\phi(t_1) & = 0 \\
\phi(t_2) & = v.
\end{align*}
Then $y(t) \neq 0$, and therefore the tangent vector to $\G$ at $x$
corresponding to $\phi$ is not in $H$.  We conclude that $H$ is not
tangent to $\G$ at $x$.
\qed
\end{Proof}

\paragraph
From now on assume that $n = \dim V$ is odd.  The Pfaffian locus $\Pf
\subset \pj^*$ is defined to be the locus of non-zero degenerate
two-forms on $V$.  (A form is called degenerate if its rank is less
than $n-1$.)  If we let $A$ be a generic $n\times n$ matrix of linear
forms on $V$, then $\Pf$ is cut out by the $n$ Pfaffians of $(n-1)
\times (n-1)$ diagonal minors of $A$ (obtained by removing the $i$-th
row and column of $A$, for $i=1,\ldots, n$).

Since a point $y\in \Pf$ corresponds to a two-form on $V$, its kernel is
a subspace $K\subset V$.  If $V$ is odd-dimensional, $y$ being
degenerate implies $K$ is at least three-dimensional.  Since we've
eliminated the zero form by projectivizing, if $\dim V = 7$, $\dim K$
could only be three or five (the rank of $y$ can be either four or two).

\paragraph
Proposition~\ref{prop:hyptgrass} shows that the Pfaffian is precisely
the classical dual variety to $\G$: indeed, $H\cap \G$ is singular for
a hyperplane $H$ if and only if $H$ is tangent at some point of $\G$,
i.e., as a point of $\pj^*$, it corresponds to a form with kernel of
dimension $\geq 3$.

\paragraph
Unlike the Grassmannian, which is smooth, the Pfaffian is singular at
the points where the rank drops further (e.g., for $\dim V = 7$, when
the rank is two). At singular points $y\in \Pf$, we have $\T_{\Pf, y} =
\T_{\pj^*, y}$.  If $y$ is smooth, the following proposition describes
the tangent space $\T_{\Pf, y} \subset\T_{\pj^*, y}$.

\begin{Proposition}
\label{prop:tpfaff}
Let $y\in \Pf$ correspond to a degenerate two-form on $V$ with kernel
$K\subset V$ of dimension three.  Identify, up to scalars, $T_{\pj^*, y}$
with $\wedge^2 V^* / \langle y \rangle$, so that vectors in $T_{\pj^*,
  y}$ are thought of as two-forms modulo $y$.  Then $v\in T_{\pj^*,
  y}$ is tangent to $\Pf$ at $y$ if and only if $K$ is isotropic for
$v$, i.e.,
\[ v(k_1\wedge k_2) = 0 \mbox{ for all } k_1, k_2 \in K. \]
\end{Proposition}

\begin{Proof}
  Pick a nonzero element $y_0\in \wedge^2 V^*$ in the cone over the
  Pfaffian that maps to $y$ under projectivization.  Since all rank
  $n-3$ skew forms on $V$ are in the same orbit of the $\GL(V)$
  action, $y_0$ is contained in the dense open subset of the cone over
  the Pfaffian which is its orbit under the $\GL(V)$ action. As a
  result, the tangent space to $\Cone(\Pf)$ at $y_0$ can be viewed as
  the span of $gy_0$ for $g\in \Lie(\GL(V))=\End(V^*)$ under the
  natural action. In an appropriate basis of $V^*$ one has
  $y_0=x_1\wedge x_2+\ldots +x_{n-4}\wedge x_{n-3}$, and
  $gy_0=gx_1\wedge x_2 + x_1\wedge gx_2 + \ldots +gx_{n-4}\wedge
  x_{n-3} + x_{n-4}\wedge gx_{n-3}$.  Since $x_i\in \Ann(K)$, we see,
  after passing from $\Cone(\Pf)$ to $\Pf$, that the tangent space to
  $\Pf$ is contained in the space of forms that make $K$ isotropic.
  It is easy to see that the Pfaffian and this space are of the same
  dimension, which finishes the proof.  \qed
\end{Proof}

\section{Simultaneous smoothness}
\label{sec:smoothness}

From now on, the space $V$ is assumed to be seven-dimensional.  In
this section we argue that, for a given choice of seven-dimensional
linear subspace
\[ W \subset \wedge^2 V^*, \]
the linear sections 
\[ X = \PM \cap \G  \]
and 
\[ Y = \PW \cap \Pf \]
are either both smooth, or both singular, at least in the case when
$X$ and $Y$ are of the (expected) dimension three.  

\paragraph
\label{subsec:defR}
In the previous section we noted that a point $x\in X$ corresponds to
a plane $T \subset V$, while a point $y\in Y$, regarded as a
two-form on $V$, has a kernel $K\subset V$, of dimension three
or five.  

Let $R\subset X \times Y$ denote the locus of pairs $(x,y)$ for which
$T \subset K$.  Let $\pi_X, \pi_Y$ be the projections from $R$
to $X$ and to $Y$, respectively.

\begin{Proposition}
The following statements hold:
\begin{itemize}
\item[(a)] the set of points $x\in X$ where $\dim T_{X,x} > 3$
  coincides with the image $\pi_X(R)$;
\item[(b)] the set of points $y\in Y$ where $\dim T_{Y,y} > 3$
  coincides with the image $\pi_Y(R)$.
\end{itemize}
\end{Proposition}

\begin{Proof}
(a) Let $x\in \G$ be a point on the Grassmannian, and $y$ an
arbitrary point in $\pj^*$.  Let $H$ be the hyperplane in
$\pj$ corresponding to $y$.  By Proposition~\ref{prop:hyptgrass}, the
following two statements are equivalent: 
\vspace*{1mm}
\begin{itemize}
\item[--] $H$ is tangent to $\G$ at $x$; 
\item[--] $T \subset K = \Ker(y)$.
\end{itemize}
\vspace*{1mm}

\noindent
Assume first that $x$ is a point in $X$ such that $\dim T_{X, x} > 3$.
Then $T_{X,x}$, which is the intersection of $T_{\G, x}$ and
$\PM$ inside $T_{\pj, x}$, has dimension higher than the
expected
\[ 3 = 10+13-20. \] 
Therefore $T_{\G, x}$ and $\PM$ do not span all of $T_{\pj, x}$, and
there exists a hyperplane $H$ inside $T_{\pj, x}$ (and thus a
hyperplane in $\pj$ through $x$) containing both.  Let $y$ be the
point in $\pj^*$ corresponding to $H$.

Since $H$ is tangent to $\G$ at $x$, $T \subset K$ by the
claim.  Therefore $y\in \Pf$ (because it has a kernel of dimension
$\geq 2$).  Since we also have $\PM \subset H$, it follows that
$y\in \PW$, and therefore $y\in Y$.  The pair $(x,y)$ is thus in
$R$, and $x$ is in the image $\pi_X(R)$.

Conversely, assume that $x$ in $X$ is in the image of $\pi_X$, and
let $(x,y)$ be a point in $R$.  Let $H$ be the hyperplane in $\pj$
corresponding to $y$.  It then follows that $H$ contains $\PM$
(because $y\in \PW$), and therefore $x\in H$.  The claim implies now
that $H$ is also tangent to $\G$ at $x$, and therefore the tangent
space to $H$ at $x$ contains $T_{\G, x}$ and $T_{\PM, x}$.  Thus
the intersection of these two spaces (which is $T_{X,x}$) can not be
of the expected dimension, and therefore $\dim T_{X,x} > 3$. 

(b) Let $y$ be a point on $Y$.  If the rank of $y$, regarded as a
two-form on $V$, is two, then $\dim T_{Y,y}>3$ because the Pfaffian is
already singular at $y$, and cutting it down by a codimension $14$
linear space $\PW$ will not cut down the dimension of the tangent
space to three.  Thus we need to show that $y$ is in the image of $R$, in
other words, that there exists a two dimensional space $T
\subset K$ such that $x\in X$.  The kernel $K$ is five-dimensional,
and thus $\G(2, K)$ is a subvariety of $\G$ of dimension six,
completely contained in the hyperplane in $\pj$ corresponding to
$y$.  Completing $\{y\}$ to a basis of $W$ gives rise to six more
hyperplanes in $\pj$, which must have a common point $x$ in
$\G(2, K)$.  Thus $x$ is in $X$ (being on $\G$ and at the
intersection of the seven hyperplanes corresponding to $W$), and
$T \subset K$.  It follows that $(x,y)\in R$, thus $y$ is in
the image of $\pi_Y$.

We can assume thus that $y$ is a smooth point of $\Pf$, the rank of
$y$ is four, and thus its kernel $K\subset V$ has dimension three.  The
tangent space $T_{\Pf, y}$ consists of all the tangent vectors in
$T_{\pj^*, y}$ for which $K$ is isotropic
(Proposition~\ref{prop:tpfaff}).  Explicitly, we have
\[ T_{\pj^*, y} \iso (\wedge^2 V^*)/\langle y\rangle, \]
and a two-form $\omega\in \wedge^2 V^*$ will be tangent to $\Pf$ at
$y$ if and only if it vanishes on $\wedge^2 K$. In other words, it
is (modulo $y$) the subspace in $\wedge^2 V^*$ whose annihilator is
$\wedge^2 K$.  

The statement $\dim T_{Y, y} > 3$ is equivalent to $T_{\PW, y}$ and
$T_{\Pf, y}$ not intersecting transversely in $T_{\pj^*, y}$,
i.e., not spanning the full $T_{\pj^*, y}$.  Since we have
\[ \Ann(T_{\PW, y}) \cap \Ann(T_{\Pf, y}) = \Ann(\langle T_{\PW,
  y}, T_{\Pf, y}\rangle), \]
it follows that this is equivalent to the existence of a non-zero
\[ x\in \Ann(T_{\PW, y}) \cap \Ann(T_{\Pf, y}) \cap \Ann(y). \]
Being in $\Ann(T_{\PW, y})\cap \Ann(y)$ is equivalent to being in
$\PM$, while being in $\Ann(T_{\Pf, y})$ is equivalent to being a
point on $\G$ for which $T \subset K$.  Thus $\dim T_{Y,y} >
3$ is equivalent to the existence of an $x\in X$ such that $(x,y) \in
R$. 
\qed
\end{Proof}

\begin{Corollary}
\label{cor:simsm}
Assume either $X$ or $Y$ has dimension three.  Then $X$ is smooth if and
only if $Y$ is smooth.
\end{Corollary}

\begin{Proof}
Both statements are equivalent to $R = \emptyset$. 
\qed
\end{Proof}

\paragraph
{\bf Remark.} 
In~\cite{Rod}, R\o dland argues that $h^{1,1}(X) = h^{1,1}(Y) = 1$ and
$h^{1,2}(X) = h^{1,2}(Y) = 50$ for generic cuts.  In fact, this
statement holds whenever $X$ and $Y$ are smooth: indeed, the family of
such cuts is smooth over the appropriate open subset of $\G(7,
\wedge^2 V^*)$, and thus all fibers are diffeomorphic.

\section{A vanishing result for Schubert cycles}
\label{sec:gvan}

In this section we define a family $\{S_y\}$ of Schubert cycles in
$\G$, parametrized by the smooth points $y$ of $\Pf$, and we study
homological orthogonality properties of this family.  The most
important result is $\G$-vanishing, which is the statement that for
$0\leq j \leq 5$, and $y_1, y_2 \in \Pf^\smooth$, we have
\[ \Ext^{j+2}(\cO_{S_{y_1}}, \cO_{S_{y_2}}(-j-1)) = 0. \]

\paragraph
Let $K\subset V$ be a linear subspace.  Define $S\subset \G$
to be the locus of two-planes $T\subset V$ (i.e., points $T\in \G$) which
intersect $K$ non-trivially.  Note that this is precisely the
Schubert cycle corresponding to the increasing sequence $(0, \dim K,
7)$.  The obvious $\GL(V)$ action permutes the Schubert cycles
corresponding to various subspaces $K$ of $V$.

\begin{Proposition}
\label{prop:defS}
  Let $K\subset V$ be a linear subspace.  Regard $\wedge^2 \Ann K
  \subset \wedge^2 V^*$ as a set of linear equations on $\pj$.  Then
  the set of two-planes $T\subset V$ such that $T\cap K \neq 0$ is
  precisely the set of closed points in $\G$ cut out by $\wedge^2 \Ann
  K$.
\end{Proposition}

\begin{Proof}
We have the following sequence of equivalent statements:
\begin{itemize}
\item[--] $T \cap K \neq 0$;
\item[--] the image $\overline{T}$ of $T$ in $V/K$ has dimension at most one;
\item[--] $\wedge^2 \overline T = 0$ in $\wedge^2 V/K$;
\item[--] for every $\overline w\in \wedge^2(V/K)^*$ we have $\overline 
  w(\wedge^2 \overline T) = 0$;
\item[--] for every $w\in \wedge^2 \Ann K$ we have $w(\wedge^2 T) = 0$.
\end{itemize}
(Note that the image of the map $\wedge^2(V/K)^* \ra \wedge^2 V^*$ is
precisely $\wedge^2 \Ann K$.)
\qed
\end{Proof}

\paragraph
The above proposition shows that $S$ is a closed subset of $\G$, and
we endow it with the reduced induced scheme structure.

\paragraph
From now on, assume that $\dim K = 3$.  Because of the $\GL(V)$
action, we can assume that in a basis $e_1,\ldots,e_7$ of $V$, $K$ is
spanned by $e_5, e_6, e_7$.  Let $x_1, \ldots, x_7$ be the dual basis of
$V^*$.  The subspace $\Ann K \subset V^*$, is spanned by $x_1, \ldots, x_4$,
and thus $\wedge^2 \Ann K$ has basis
\[ x_1\wedge x_2, x_1\wedge x_3, \ldots, x_3\wedge x_4. \]

\begin{Proposition}
\label{prop:cutS}
The Schubert cycle $S$ is cut out scheme-theoretically by $\wedge^2
\Ann K$.  It is a rational variety with rational Cohen-Macaulay
singularities, of codimension three in $\G$.
\end{Proposition}

\begin{Proof}
In the Zariski open subspace $U\iso \C^{10}\subseteq \G$ in which
points $T$ are given by
\[ T=\Span(a_{11}e_1+\ldots+a_{15}e_5+e_6,
a_{21}e_1+\ldots+a_{25}e_5+e_7), \] 
the cycle $S$ is characterized by the condition that the matrix
$R=\{a_{ij}\}$, $i=1,2$, $j=1,\ldots,4$ has rank one.  The equations
from $\wedge^2 \Ann K$ are spanned by $x_i\wedge x_j$, $1\leq i<j \leq
4$, which are precisely the maximal minors of the matrix $R$.  It is a
well-known fact that this determinantal variety is reduced and
irreducible, and it coincides with the product of $\C^2$ and the cone
over the Segre embedding of $\pj^1\times \pj^3$.  Thus $S \cap U$ is
cut out by $\wedge^2 \Ann K$ and is rational.  It remains to observe
that every point in $\G$ can be mapped inside $U$ by an element of
$\GL(V)$ that fixes $K$.

The singularities of $S$ are rational Cohen-Macaulay since $S\cap U$
is a toric variety (also see Propositions~\ref{prop:resS}
and~\ref{prop:resOS}).   
\qed
\end{Proof}

\begin{Proposition}
\label{prop:resOS}
Denote by $T$ the rank two tautological bundle on $\G$, and let $A =
\Ann K$.  Then there exists a resolution of $\cO_S$ on $\G$ of the form:
\[ 0 \ra \wedge^4 A \otimes \Sym^2(T)(-1) \ra \wedge^3 A \otimes T(-1)
\ra \wedge^2 A \otimes \cO_\G(-1) \ra \cO_\G \ra \cO_S \ra 0. \]
\end{Proposition}

\begin{Proof}
In the local chart described in the proof of
Proposition~\ref{prop:cutS}, the above resolution becomes the
Eagon-Northcott complex resolving the locus of rank one matrices
inside the space of $2\times 4$ matrices
\[ \left ( \begin{array}{llll} a_{11} & a_{12} & a_{13} & a_{14} \\
a_{21} & a_{22} & a_{23} & a_{24} \end{array} \right ). \]
\qed
\end{Proof}

\paragraph
We are now interested in getting a better understanding of the
geometry of $S$.  Let $\tS$ be the set of pairs $(T_1,T_2)$ of linear
subspaces of $V$, of dimensions $1$ and $2$, respectively, such that
$T_1 \subseteq T_2 \cap K$.  Since $T_1 \subset K$ is a
one-dimensional linear subset, the first projection $\mu$ maps $\tS$
to $\pj(K)$.  Similarly, the second projection $\pi$ maps $\tS$ to
$S$.  Thus we have the diagram
\[
\begin{diagram}[height=2em,width=2em]
\tS & \rTo^\mu & \pj(K) \\
\dTo^\pi & & \\
S.
\end{diagram}
\]
Giving a point $(T_1, T_2) \in \tS$ amounts to choosing a line $T_1
\subset K$, along with another line $T_2/T_1$ in the six-dimensional
space $V/T_1$.  Expressing this globally, we see that $\tS$ is the
projectivization of the rank six vector bundle on $\pj(K)$
\[ \big ( V \otimes \cO \big ) / \cO(-1). \]

\begin{Proposition}
\label{prop:resS}
The map $\pi$ is a resolution of singularities of $S$.  There is a
short exact sequence of bundles on $\tS$
\[ 0 \ra \mu^*\cO_{\pj(K)}(-1) \ra \pi^* T \ra \cO_{\mu}(-1) \ra 0, \]
where $T$ is the tautological rank two vector bundle on $\G$,
restricted to $S$, and $\cO_\mu(-1)$ is the tautological relative
bundle of the projection $\mu$.  We have
\[ \pi^* \cO_S(-1) = \mu^* \cO_{\pj(K)}(-1) \otimes \cO_\mu(-1). \]
\end{Proposition}

\begin{Proof}
The space $\tS$ is smooth, and $\pi$ is obviously birational, thus
$\pi$ is a resolution of singularities.  The short exact sequence is
nothing but a global way of expressing the fact that $T_2$ is an
extension of $T_1$ and $T_2/T_1$.  The last equality follows from the
fact that on $\G$ we have
\[ \wedge^2 T = \cO_\G(-1). \]
\qed
\end{Proof}

\begin{Proposition}
\label{prop:vanS}
For $3\leq j \leq 5$, $0\leq k \leq 2$, the bundles $\cO(j-6)$ and
$\Sym^k(T^*)(-j)$ have no cohomology on $S$. 
\end{Proposition}

\begin{Proof}
Since the singularities of $S$ are rational, we can work on $\tS$ by
replacing all the bundles by their pull-backs by $\pi$.  Using
Proposition~\ref{prop:resS} we see that that $T^*$ is an extension of
$\mu^* \cO_{\pj(K)}(1)$ and $\cO_\mu(1)$.  Thus $\Sym^2(T^*)$ has a
filtration with associated graded pieces
\[ \mu^* \cO_{\pj(K)}(2), \mu^* \cO_{\pj(K)}(1) \otimes \cO_\mu(1),
\cO_\mu(2). \] 
Using the computation of $\cO(-1)$ from the same Proposition, it
follows that all the bundles we are interested in have filtrations
whose associated graded object is a sum of line bundles of the form
\[ \mu^* \cO_{\pj(K)}(k) \otimes \cO_\mu(l), \] with $-6\leq l\leq
-1$.  By the Leray spectral sequence, these line bundles have no
cohomology.  
\qed
\end{Proof}

\paragraph
\label{subsec:defS}
We are now ready to define the Schubert cycles $S_y$, parametrized by
the smooth points of $\Pf$.  Let $y \in \Pf$ be a smooth point, which
is thought of as a two-form on $V$ of rank four.  Let $K_y$ be the
kernel of this form, a linear three-space in $V$.  We define the locus
$S_y \subset \G$ to be the Schubert cycle $S$ defined above, for
$K=K_y$.

\begin{Proposition}
\label{prop:globalS}
There exists a subvariety $\mathbf{S}$ of $\G \times \Pf^\smooth$,
flat over $\Pf^\smooth$, whose fiber $\mathbf{S}_y\subset \G$ over any
closed point $y\in\Pf^\smooth$ is precisely $S_y$.
\end{Proposition}

\begin{Proof}
Fix a point $y_0\in \Pf^\smooth$, and let $S_0$ be the corresponding
subvariety of $\G$.  The group $\GL(V)$ acts transitively on both $\G$
and $\Pf^\smooth$.  One can take $\mathbf{S}$ to be the image of
$\GL(V) \times S_0$ under this action, with the reduced-induced scheme
structure.  The flatness statement follows from~\cite[Theorem
9.9]{HarAG}, as all translates of $S_0$ have the same Hilbert
polynomial in $\G$.  
\qed
\end{Proof}

\begin{Proposition}{\bf ($\G$-vanishing).}
\label{prop:gvan}
Let $S_1$, $S_2$ be the Schubert varieties defined
in~(\ref{subsec:defS}), corresponding to linear three-dimensional
subspaces $K_1$, $K_2$ of $V$, and let $\cI_{S_1}$, $\cI_{S_2}$ denote
their ideal sheaves in $\G$.  Then, for $0\leq j\leq 5$, we have
\[ \Ext^{\bullet}_{\G}(\cI_{S_1}, \cI_{S_2}(-j-1)) = 0. \]
\end{Proposition}

\begin{Proof}
In view of Serre duality, it suffices to show that for $j=3$, $4$, or
$5$, we have
\[ \Ext^\bullet_{\G}(\cI_{S_1}, \cI_{S_2}(-j-1)) = 0. \] 
From the long exact sequence
\[ 0 \ra \cI_{S_2}(-j-1) \ra \cO_\G(-j-1) \ra \cO_{S_2}(-j-1) \ra 0 \]
we will get the desired vanishing if we prove that 
\[ \Ext^\bullet_\G(\cI_{S_1}, \cO_\G(-j-1)) =
\Ext^\bullet_\G(\cI_{S_1}, \cO_{S_2}(-j-1)) = 0. \] 
By Serre duality and standard facts about the cohomology of $\cO_\G$,
the first statement amounts to proving that
\[ H^\bullet(\G, \cI_{S_1}(j - 6)) = H^\bullet(S_1, \cO_{S_1}(j-6)) =
0, \]
which is part of Proposition~\ref{prop:vanS}.  

On the other hand, the resolution from Proposition~\ref{prop:resOS}
together with the local-to-global spectral sequence will give us the
second vanishing statement provided we argue that 
\[ H^\bullet(S_2, \Sym^k(T^*)(-j)) = 0 \]
for $0\leq k \leq 2$.  This is also part of
Proposition~\ref{prop:vanS}.
\qed
\end{Proof}

\paragraph
\label{subsec:Kuzvan}
As pointed out by Kuznetsov, $\G$-vanishing also follows
from~\cite[Theorem 3.1]{Kuzsod}, using the resolution obtained in
Proposition~\ref{prop:resOS}.

\section{A family of curves}
\label{sec:allcurves}

In this section we define, for a point $y\in Y$, a curve $C_y$ in $X$.
(We abuse the notation slightly: $C_y$ may not be reduced or
irreducible, but it does have dimension one.)  The family
$\{\cI_{C_y}\}_{y\in Y}$ of ideal sheaves of these curves is the
orthogonal family which induces the equivalence of derived categories
$\D(X) \iso \D(Y)$.  We then argue that $\dim C_y = 1$ for every
choice of $y\in Y$ (which is essential in proving that the family of
$C_y$'s is flat), and that
\[ y_1 \neq y_2 \implies \Hom_X(\cI_{C_{y_1}}, \cI_{C_{y_2}}) = 0, \]
which is later needed for the orthogonality of the family.

Most of our statements depend on the assumption that $X$ and $Y$ are
smooth.  Therefore, for the remainder of this paper we shall assume
that a space $W \subset \wedge^2 V^*$ is chosen in such a way that $Y$
(and therefore $X$) is smooth.

\begin{Lemma}
\label{lem:diffker}
Let $W\subset \wedge^2 V^*$ be a seven-subspace, and assume that $Y = \PW
\cap \Pf$ is smooth of dimension three.  Then the kernels of any pair of points $w_1 \neq w_2$ in $Y$ (thought of as two-forms on $V$) are distinct.
\end{Lemma}

\begin{Proof}
  Assume that $\Ker w_1 = \Ker w_2 = K$.  All the linear combinations
  of $w_1$ and $w_2$ have $K$ contained in their kernel, so we can
  think of them as a one-dimensional family of forms on the
  four-dimensional space $V/K$.  At some point $w$ in this family the
  rank will drop to two, as this is a codimension one condition.  At
  $w$ the space $\Pf$ is singular, and therefore the regular cut $Y$
  is singular at $w$ as well.  \qed
\end{Proof}

\begin{Lemma}
\label{lem:kerw}
Assume that $X$ is smooth of dimension three.  Then the space of
linear forms $w\in \pj^*$ which vanish along all of $X$ is precisely
$W$.
\end{Lemma}

\begin{Proof}
Another way to phrase this assertion is by saying that $W$ is the
kernel of the map
\[ \rho: H^0(\pj, \cO_\pj(1)) \ra H^0(X, \cO_X(1)). \] Since the space
$W$ is already contained in the kernel of $\rho$, it suffices to show
that $\dim \Ker \rho = 7$.  Factor $\rho$ into the composition
\[ H^0(\pj, \cO_\pj(1)) \ra H^0(\G, \cO_\G(1)) \ra H^0(X,
\cO_X(1)), \] 
where it is well-known that the first map is an isomorphism.

Consider the ideal sheaf $\cI_X$ of $X$ in $\G$.  Its twist by one has
the Koszul resolution
\[ 0 \ra \cO_\G(-6) \ra \cdots \ra \cO_\G(-1)^{21} \ra \cO_\G^7 \ra
\cI_X(1) \ra 0, \]
where the sheaves $\cO_\G(-1), \ldots, \cO_\G(-6)$ are acyclic with no
global sections.  It follows that 
\[ H^0(\G, \cI_X(1)) = H^0(\G, \cO_\G^7) = \C^7. \]
From this and the short exact sequence 
\[ 0 \ra \cI_X(1) \ra \cO_\G(1) \ra \cO_X(1) \ra 0 \]
it follows that $\dim\Ker \rho = 7$.
\qed
\end{Proof}

\paragraph
\label{subsec:defC}
Let $\mathbf{C}$ be the scheme theoretic intersection of $\mathbf{S}$
with $X\times Y$ inside $\G\times \Pf$, regarded as a subscheme of
$X\times Y$.  The sheaf $\cI_{\mathbf{C}}$, the ideal sheaf of
$\mathbf{C}$ in $X\times Y$, regarded as an object of $\D(X\times Y)$,
will be the kernel of the Fourier-Mukai transform $\D(Y) \ra \D(X)$.

\begin{Proposition}
\label{prop:allcurves}
  The fiber of $\mathbf{C}$ over any closed point $y\in Y$ has
  dimension one.  The scheme $\mathbf{C}$ is flat over $Y$.
\end{Proposition}

\begin{Proof}
By an appropriate version of Bertini's theorem, we can choose a
sequence of cuts going from $\G$ to $X$ so that all the intermediate
members are smooth.  Lefschetz's theorem then gives $\Pic X = \Z$,
generated by the restriction of $\cO_\G(1)$.  For a point $y\in Y$,
$C_y$ is cut out from $X$ by the same six linear equations that cut
out $S_y$ from $\G$, i.e., $C_y$ is the result of intersecting $X$
with the hyperplanes corresponding to the points of the
six-dimensional linear space $L = \wedge^2 \Ann \Ker y$.  Every form
$w\in L$ contains $\Ker y$ in its kernel.

Assume by contradiction that $C_y$ contains a divisor $D$.  We claim
that the intersection $H\cap X$ is either $D$, or all of $X$, for
every hyperplane $H$ in $L$.  Indeed, any non-trivial intersection
$H\cap X$ is the generator of $\Pic X$, and as such it is a prime
divisor on $X$.  Therefore any such cut must equal $D$, since it
already contains it.

Pick a point $P$ in $X\setminus D$.  The condition that $H$ contain
$P$ is a single linear condition on $L$, thus there is a five-dimensional
space $R$ of forms in $L$ that vanish at $P$.  The above discussion
shows that if a form in $L$ vanishes at a point outside $D$, then it
vanishes on all of $X$.  Thus the forms in $R$ vanish on all of
$X$. 

By Lemma~\ref{lem:kerw}, then, we have $R\subset W$.  Any two
linearly-independent forms of rank four in $R$ yield two points in $Y$
whose kernel (as forms) is $K$.  By Lemma~\ref{lem:diffker}, this
contradicts the assumption that $Y$ is smooth.  Therefore our
assumption that $C_y$ contains a divisor must be false.

The subvarieties $S_y\subset \G$ for $y\in Y$ are all isomorphic, and
they have the same Hilbert polynomial (in fact, there is a transitive
$\GL(V)$-action permuting them).  If we let $y, y_1,\ldots,y_6$ be a
basis of $W$, them $C_y$ is obtained from $S_y$ by six cuts with
hyperplanes corresponding to $y_1, \ldots , y_6$.  Since $C_y$ is a
curve (codimension six in $W$), it follows that this is a regular
sequence of cuts, and thus the Hilbert polynomial of $C_y$ does not
depend on $y\in Y$.  Therefore the family $\{C_y\}_{y\in Y}$ is flat
over $Y$, by~\cite[III.9.9]{HarAG}.
\qed
\end{Proof}

\begin{Proposition}
\label{prop:Hom0}
Let $y$ be a point in $Y$.  Then we have
\[ \Hom_X(\cI_{C_{y}}, \cI_{C_{y}}) = \C. \]
\end{Proposition}

\begin{Proof}
From the short exact sequence 
\[ 0 \ra \cI_{C_y} \ra \cO_X \ra \cO_{C_y} \ra 0 \]
we get 
\[ 0 \ra \C \ra \Hom_X(\cI_{C_y}, \cO_X) \ra \Ext^1_X(\cO_{C_y},
\cO_X) = H^2(X, \cO_{C_y})^* = 0. \]
Thus any non-zero homomorphism $\cI_{C_y} \ra \cO_X$ is, up to a
scalar multiple, the usual inclusion.

Let $f:\cI_{C_y} \ra \cI_{C_y}$ be any non-zero homomorphism.
Composing with the inclusion $\cI_{C_y} \hookrightarrow \cO_X$ we get
a non-zero map $\cI_{C_y} \ra \cO_X$ which, by the above reasoning,
must be the usual inclusion of $\cI_{C_y}$ into $\cO_X$.  Therefore
$f$ is injective, and since $\cI_{C_y}$ is torsion-free, it follows that
$f$ is a multiple of the identity map.
\qed
\end{Proof}
\section{A reduction to the Grassmannian}
\label{sec:redGrass}

In this section we argue that checking the orthogonality of the family
$\{\cI_{C_y}\}_{y\in Y}$ can be reduced to the $\G$-vanishing
statement (Proposition~\ref{prop:gvan}).

\paragraph 
We assume that we have chosen $W$ so that $X$ and $Y$ are smooth of
dimension three.  Let $y_1$ and $y_2$ be two distinct points in $Y$.
Regarded as smooth points of the Pfaffian $\Pf$, $y_1$ and $y_2$ give
rise to Schubert cycles $S_1$ and $S_2$ in $\G$ as explained in
(\ref{subsec:defS}).  These are rational, integral, codimension three
subschemes of $\G$, whose singularities are rational and
Cohen-Macaulay.

Let $C_1$ and $C_2$ be the corresponding curves in $X$ obtained by
taking scheme-theoretic intersections of $S_1$ and $S_2$ with $X$,
respectively, as defined in Section~\ref{sec:allcurves}.

We recall the $\G$-vanishing statement of Proposition~\ref{prop:gvan}:
\[ \Ext^\bullet_\G(\cI_{S_1}, \cI_{S_2}(-j-1)) = 0 \mbox { for } 0 \leq
j \leq 5. \] 

\paragraph
\textbf{Notation.}  
Throughout this section, if $S$ is a subscheme of a scheme $Z$, we
shall denote by $\cI_S$ the ideal sheaf of $S$, regarded as a coherent
sheaf on $Z$.  We shall always make sure to be precise, if $Z$ is
itself a subscheme of another scheme $Z'$, which scheme is $S$
regarded as a subscheme of, $Z$ or $Z'$.

\begin{Proposition}
\label{prop:redGrass}
For $y_1 \neq y_2$ we have
\[ \Ext^\bullet_X(\cI_{C_1}, \cI_{C_2}) = 0. \]
\end{Proposition}

\paragraph
Before we begin the proof of Proposition~\ref{prop:redGrass} we need
several intermediate vanishing results.

\begin{Lemma}
\label{lem:van1}
For $0\leq j \leq 5$ we have
\[ H^{\bullet}(\G, \cI_{S_2}(-j)) = 0. \]
\end{Lemma}

\begin{Proof}
Consider the short exact sequence
\[ 0 \ra \cI_{S_2} \ra \cO_{\G} \ra \cO_{S_2} \ra 0. \]
The cohomology groups
\[ H^{\bullet}(\G, \cO_{\G}(-j)) \] 
vanish for $1\leq j \leq 6$, by Kodaira vanishing.  Since $S_2$ is
irreducible, rational, and with rational singularities, we conclude that
\[ H^\bullet(\G, \cI_{S_2}) = 0, \]
which is the $j=0$ case.  For $1\leq j \leq 5$ we have
\[ H^{\bullet+1}(\G, \cI_{S_2}(-j)) = H^\bullet(S_2, \cO_{S_2}(-j)), \]
which vanishes by Proposition~\ref{prop:vanS}.
\qed
\end{Proof}

\paragraph
We now divide the proof of Proposition~\ref{prop:redGrass} into
several steps, corresponding to decreasing subvarieties of $\G$
obtained by successive cuts with linear subspaces of increasing
codimension, corresponding to subspaces of $\PW$ of increasing
dimension.

To begin, let $H_1$ be the hyperplane in $\pj$ corresponding to $y_1$,
and let $Z_1$ be the scheme-theoretic intersection of $\G$ and
$H_1$. It is a hypersurface of $\G$ whose singularities are locally
isomorphic to the product of $\C^2$ and an ordinary double point of
dimension seven.  Consider $L_1$ and $L_2$ which are the intersection
of $S_1$ and $S_2$ with $H_1$, regarded as subschemes of $Z_1$.  While
$H_1$ cuts down the dimension of $S_2$ by one, it already contains
$S_1$.  (The fact that $H_1$ cuts down the dimension of $S_2$ follows
from Proposition~\ref{prop:allcurves}.)  Thus $L_1$ and $L_2$ have
codimensions two and three, respectively, in $Z_1$.

\begin{Proposition}
\label{prop:Z1van}
$\G$-vanishing implies $Z_1$-vanishing:
\[ \Ext^{\bullet}_{Z_1}(\cI_{L_1}, \cI_{L_2}(-j)) = 0 \mbox{ for } 0 \leq
j \leq 5. \]
\end{Proposition}

\begin{Proof}
Let $r:Z_1\ra \G$ denote the natural embedding.  Since $L_2$ is
obtained by a transversal cut of $S_2$ by $H_1$, we have $\cI_{L_2} =
r^* \cI_{S_2}$.  By an easy form of Grothendieck duality it follows that
\begin{align*}
\Ext^{\bullet}_{Z_1}(\cI_{L_1}, \cI_{L_2}(-j)) & =
\Ext^{\bullet}_{Z_1}(\cI_{L_1}, r^* \cI_{S_2}(-j)) \\
& = \Ext^{\bullet}_{\G}(r_!\cI_{L_1}, \cI_{S_2}(-j)) \\
& = \Ext^{\bullet+1}_{\G}(r_*\cI_{L_1}, \cI_{S_2}(-j-1)),
\end{align*}
since the embedding $Z_1\subset \G$ has relative dualizing complex
$\cO_{Z_1}(1)[-1]$.

Consider the short exact sequence of sheaves on $\G$
\[ 0 \ra \cI_{Z_1}=\cO_\G(-1) \ra \cI_{S_1} \ra r_* \cI_{L_1} \ra 0. \] 
Writing down the long exact sequence of $\Ext$'s we get the result
from $\G$-vanishing and the vanishing of
\[ \Ext^{\bullet}_{\G}(\cI_{Z_1}, \cI_{S_2}(-j-1)) = H^{\bullet}(\G,
\cI_{S_2}(-j)), \]
which is Lemma~\ref{lem:van1}.
\qed
\end{Proof}

\begin{Lemma}
\label{lem:van2}
For $0\leq j \leq 5$ we have
\[ H^{\bullet}(Z_1, \cI_{L_1}(-j)) = 0. \]
\end{Lemma}

\begin{Proof}
  The proof is essentially the same as that of Lemma~\ref{lem:van1}
  and will be left to the reader.  One needs to use the fact that
  $Z_1$ is rational with rational singularities and $K_{Z_1} = \cO(-6)$. 
  \qed
\end{Proof}

\paragraph
We can now complete the proof of Proposition~\ref{prop:redGrass}.  Let
$H_2$ be the hyperplane in $\pj$ corresponding to $y_2$, and let
$Z_{1,2} = \G \cap H_1 \cap H_2$.  Let $D_1$, $D_2$ be the
intersections of $S_1$ and $S_2$ with $Z_{1,2}$, regarded as
subschemes of $Z_{1,2}$.  Note that since $S_2 \subset H_2$, we have
$D_2 = L_2$; both $D_1$ and $D_2$ are now codimension two in $Z_{1,2}$.

Let $W_5$ be any linear subspace of $W$ which, together with the
one-dimensional subspaces corresponding to $y_1$ and $y_2$, spans $W$.
Because of dimension reasons, $X$, $C_1$, and $C_2$ are obtained from
$Z_{1,2}$, $D_1$, and $D_2$ by five successive transversal cuts with
five hyperplanes corresponding to a basis of $W_5$.  

The strategy of the proof is to look at the successive embeddings 
\[ X \stackrel{g}{\lra} Z_{1,2} \stackrel{h}{\lra} Z_1
\stackrel{r}{\lra} \G, \] 
and to use the vanishing of $\Ext$ groups on each one to conclude the
vanishing of other $\Ext$ groups on the previous one.  We have already
argued (Proposition~\ref{prop:Z1van}) that $\G$-vanishing implies
$Z_1$-vanishing.  The passage from $Z_1$-vanishing to the required
vanishing on $X$ will be done by an appropriate resolution of
$h_* g_* \cI_{C_2}$ on $Z_1$.

We have
\[ \cI_{C_1} = g^* \cI_{D_1}, \]
thus
\[ \Ext^\bullet_X(\cI_{C_1}, \cI_{C_2}) = \Ext^\bullet_X(g^* \cI_{D_1}, \cI_{C_2})
= \Ext^\bullet_{Z_{1,2}}(\cI_{D_1}, g_* \cI_{C_2}), \]
and 
\[ g_* \cI_{C_2} = \cI_{D_2} \otimes_{Z_{1,2}} \cO_X. \]
Because the intersection of $D_2$ and $X$ inside $Z_{1,2}$ is
transversal, there are no higher $\Tor$'s in the above tensor product.
We conclude that $g_* \cI_{C_2}$ is quasi-isomorphic to the complex
\begin{align*} 
  \cI_{D_2} \otimes & \Koszul^{Z_{1,2}}(W_5) = \\
  & (0 \ra \cI_{D_2}(-5) \otimes \wedge^5 W_5 \ra \cdots \ra \cI_{D_2}(-1)
  \otimes W_5 \ra \cI_{D_2}\ra 0).
\end{align*}

Now let us move on to the embedding $h:Z_{1,2} \ra Z_1$.  Again,
\[ \Ext^\bullet_{Z_{1,2}}(\cI_{D_1}, g_* \cI_{C_2}) = \Ext^\bullet_{Z_{1,2}}(h^*
\cI_{L_1}, g_* \cI_{C_2}) = \Ext^\bullet_{Z_1}(\cI_{L_1}, h_* g_*
\cI_{C_2}). \]
We use the previous resolution of $g_* \cI_{C_2}$ to compute a
resolution of $h_* g_* \cI_{C_2}$ on $Z_1$.  The main difference
is that $\cI_{D_2}$ is not obtained from $\cI_{L_2}$ by a transversal
intersection.  In fact, $L_2$ is already contained in $H_2$, so we
have a short exact sequence on $Z_1$
\[ 0 \ra \cI_{Z_{1,2}} = \cO_{Z_1}(-1) \ra \cI_{L_2} \ra h_* \cI_{D_1}
\ra 0, \]
and $h_* \cI_{D_1}$ is quasi-isomorphic to the two-term complex
\[ 0 \ra \cO_{Z_1}(-1) \ra \cI_{L_2} \ra 0. \]
Observe that the Koszul complex restricts from $Z_1$ to $Z_{1,2}$,
because the five cuts from $W_5$ are transversal to $Z_1$ and
$Z_{1,2}$.  In other words,
\[ h^* \Koszul^{Z_1}(W_5) = \Koszul^{Z_{1,2}}(W_5). \]
The projection formula for derived categories implies that
\begin{align*}
h_* g_* \cI_{C_2} & = h_*(\Koszul^{Z_{1,2}}(W_5) \otimes \cI_{D_1}) = 
h_*(h^*\Koszul^{Z_1}(W_5) \otimes \cI_{D_1}) \\
& = \Koszul^{Z_1}(W_5) \otimes h_* \cI_{D_1}
 = \Koszul^{Z_1}(W_5) \otimes (\cO_{Z_1}(-1) \ra \cI_{L_2}).
\end{align*}
All operations above are to be understood as derived.  The
computations of derived tensor product, and derived pull-back are
correct, as the Koszul complex is locally free.  Since $g$ and $h$ are
embeddings, the left hand side is just the one sheaf complex $h_* g_*
\cI_{C_2}$.

Thus $h_* g_* \cI_{C_2}$ is quasi-isomorphic to the total complex of 
\[ \Koszul^{Z_1}(W_5) \otimes (\cO_{Z_1}(-1) \ra \cI_{L_2}). \]
By the local-to-global spectral sequence, in order to conclude that 
\[ \Ext_X^\bullet(\cI_{C_1}, \cI_{C_2}) = 0 \]
it suffices to prove that, for $0\leq j\leq 5$, we have
\begin{align*}
& \Ext_{Z_1}^{\bullet}(\cI_{L_1}, \cI_{L_2}(-j)) = 0 
\intertext{and}
&\Ext_{Z_1}^{\bullet}(\cI_{L_1}, \cO_{Z_1}(-j-1)) = 0. 
\end{align*}
The first statement is precisely $Z_1$-vanishing, which is implied by
$\G$-vanishing by Proposition~\ref{prop:Z1van}.  The second statement
follows by Serre duality: we have
\[ \Ext_{Z_1}^{\bullet}(\cI_{L_1}, \cO_{Z_1}(-j-1)) = H^{9-\bullet}(Z_1,
\cI_{L_1}(j-5))^* \]
which vanishes for $0 \leq j \leq 5$ by Lemma~\ref{lem:van2}.\qed

\section{The derived equivalence}
\label{sec:deq}

In this section we define the integral transform $\Phi$, and we argue
that it gives an equivalence of derived categories $\D(Y) \iso \D(X)$
by verifying the criterion of~\cite[Theorem 1.1]{BriFM}. 

\paragraph
Let $\mathbf{C}$ be the dimension four subscheme of $X\times Y$ defined
in~(\ref{subsec:defC}).  We take its ideal sheaf $\cI_{\mathbf{C}}$ in
$X\times Y$ as the kernel of an integral transform
\[ \Phi: \D(Y) \ra \D(X)\quad\quad \Phi(\cE) = \R\pi_{X,*}(\pi_Y^*(\cE)
\lotimes \cI_{\mathbf{C}}). \]
(Here, $\pi_X$ and $\pi_Y$ are the projections from $X\times Y$ to $X$
and $Y$, respectively.)

\begin{Theorem}
\label{thm:mainthm1}
The functor $\Phi$ is a Fourier-Mukai transform, i.e., it is an
equivalence of categories $\D(Y) \iso \D(X)$.
\end{Theorem}

\begin{Proof}
Observe that since $\mathbf{C}$ is flat over $Y$, $\Phi\cO_y$ is
precisely the ideal sheaf $\cI_{C_y}$ on $X$.  Thus we can think of
$\{\cI_{C_y}\}_{y\in Y}$ as a family of sheaves on $X$, parametrized
by the points of $Y$.  Moreover, this family satisfies the following
properties:
\begin{itemize}
\item[1.] $\Ext^i_X(\cI_{C_{y_1}}, \cI_{C_{y_2}}) = 0$ for any $i$,
  and any pair of distinct points $y_1$, $y_2$ in $Y$
  (Proposition~\ref{prop:redGrass}); 
\item[2.] $\cI_{C_{y}}$ is a simple sheaf on $X$ for $y\in Y$
  (Proposition~\ref{prop:Hom0}); 
\item[3.] $X$ is Calabi-Yau, therefore $\cI_{C_y} \otimes \omega_X \iso
  \cI_{C_y}$.
\end{itemize}
Thus, in the terms of~\cite{BriFM}, the sheaf $\cI_{\mathbf{C}}$ is
strongly simple.  Combining this with Property 3, Theorem 1.1 of
[loc.cit.] shows that $\Phi$ is a Fourier-Mukai transform.  
\qed
\end{Proof}

\section{Connections with existing work}
\label{sec:comments}

In this section we put the example in this paper in the context of
Kuznetsov's Homological Projective Duality.

\paragraph
The derived equivalence that we obtain appears to be a particular case
of Homological Projective Duality, as explained by
Kuznetsov~\cite{Kuz}.  Let us recall the setting of this theory.

Let $X \subset \pj(V)$ be a smooth projective variety, and assume that
its derived category admits a semi-orthogonal decomposition of the form
\[ \D(X) = \langle \cA_0, \cA_1(1), \ldots, \cA_{i-1}(i-1) \rangle\]
where
\[ 0 \subseteq \cA_0\subseteq \cA_1 \subseteq \cdots \subseteq
\cA_{i-1} \] 
are full subcategories of $\D(X)$, and $(i)$ denotes twisting by
$\cO_X(i)$.  If $H$ is a hyperplane in $\pj(V)$, let $X_H$ be the
corresponding hyperplane section of $X$.  Then it is easy to see that
\[ \langle \cA_1(1), \cA_2(2), \ldots, \cA_{i-1}(i-1) \rangle \]
is a semiorthogonal collection in $\D(X_H)$.  In general there is no
reason to expect this collection to generate $\D(X_H)$.  Let $\cC_H$
denote the orthogonal in $\D(X_H)$ of the subcategory generated by the
above collection. 

We can think of $\{\cC_H~:~H \in \pj(V^*)\}$ as a family of triangulated
categories, parametrized by $H\in \pj(V^*)$.  In certain cases we can
find a smooth variety $Y$, together with a morphism $Y \ra \pj(V^*)$,
such that $\{\D(Y_H)\}_{H\in \pj(V^*)}$ is essentially the family of
$\cC_H$'s (for details, see~\cite{Kuz}).  In this situation $Y$ is
called the Homological Projective Dual of $X$.

The main theorem of~\cite{Kuz} is the statement that $Y$ then admits a
semiorthogonal decomposition of a similar form, and if $L$ is a linear
subspace of $V^*$, $L^\perp \subset V$ its annihilator, then the
linear sections 
\[ X_L = X \times_{\pj(V)} \pj(L^\perp),\quad
Y_L = Y \times_{\pj(V^*)} \pj(L) \]
have closely related derived categories $\D(X_L)$, $\D(Y_L)$.
Explicitly, their categories will have several trivial components (for
$X_L$, arising from the semiorthogonal decomposition of $X$, for $Y_L$
from that of $Y$), as well as a component which is {\em the same} in
$X_L$ and $Y_L$.  In particular, if the dimension of $L$ is chosen
properly, the trivial components will disappear, and we get 
\[ \D(X_L) \iso \D(Y_L). \]
Another important property of this setup is that $X_L$ is smooth
precisely when $Y_L$ is, and thus $Y$ is closely related to the
classical projective dual of $X$.

\paragraph
In our setting, we have the smooth projective variety $\G \subset
\pj$, its projective dual $\Pf$, and linear cuts $X$ and $Y$ of them
by dual linear spaces.  The varieties $X$ and $Y$ are simultaneously
smooth, and when they are, their derived categories are equivalent.
It is very tempting, in this context, to conjecture the following.
\begin{Conjecture}
  There exists a Lefschetz decomposition of the derived category of
  $\G$, and a smooth variety $\overline{\Pf} \ra \pj^*$ (possibly
  understood in an extended sense, for example a non-commutative
  scheme), mapping to $\Pf$, which is homologically projectively dual
  to $\G$ with respect to the decomposition of $\D(\G)$.  Furthermore,
  Theorem~\ref{thm:mainthm} is a direct application of the main result
  of~\cite{Kuz}.
\end{Conjecture}

\paragraph
The current form of the conjecture was suggested by Kuznetsov, who
pointed out that our choice of $\overline{\Pf}$ from an early version
of the paper can not be the correct homologically projective dual.  In
fact, Kuznetsov claims to have constructed such a non-commutative
homologically projective dual of $\G(2,7)$ in~\cite{Kuz1}. In
addition, he expects to be able to extend this result to $\G(2,2n+1)$
for $n>3$ as well. He hopes to address the subject further in an
upcoming paper.

\end{document}